\documentclass[11pt,a4paper,reqno]{amsart}
\usepackage[T1]{fontenc}
\usepackage{lmodern,amssymb,mathtools,microtype}
\usepackage{a4wide}
\usepackage{xcolor}
\definecolor{linkblue}{RGB}{30,65,90}
\usepackage{hyperref}
\hypersetup{colorlinks=true,linkcolor=linkblue,citecolor=linkblue,urlcolor=linkblue,
 pdftitle={A short proof of the global minimality of the cracktip for the Mumford--Shah functional},
 pdfauthor={Francesco Deangelis, Camillo De Lellis, Matteo Focardi}}
\usepackage{bookmark}
\newtheorem{theorem}{Theorem}[section]
\theoremstyle{remark}

\newcommand{\R}{\mathbb R}
\newcommand{\HH}{\mathcal H}
\newcommand{\E}{\mathcal E}
\numberwithin{equation}{section}
\title[Global minimality of the cracktip]{A remark on the global minimality of the cracktip}
\author[F. Deangelis]{Francesco Deangelis}
\address{Applied Mathematics M\"unster, University of M\"unster,
Einsteinstrasse 62, 48149 M\"unster, Germany}
\email{francesco.deangelis@uni-muenster.de}
\author[C. De Lellis]{Camillo De Lellis}
\address{School of Mathematics, Institute for Advanced Study,
1 Einstein Dr., Princeton NJ 08540, USA}
\address{Gran Sasso Science Institute, viale Francesco Crispi, 7,
67100 L'Aquila, Italy}
\email{camillo.delellis@ias.edu}
\author[M. Focardi]{Matteo Focardi}
\address{DiMaI U. Dini, Universit\`a di Firenze, Florence, Italy}
\email{matteo.focardi@unifi.it}
\date{}
\keywords{Mumford--Shah functional, cracktip, global minimizers,
Dirichlet boundary conditions, regularity}
\begin{document}
\begin{abstract}
We give a short proof that the planar cracktip is a global minimizer
of the Mumford--Shah functional, based on interior and Dirichlet
boundary regularity.
We construct a Dirichlet minimizer with nonempty discontinuity
set and show that this set contains no simple closed curves.
The regularity results then imply the existence of an interior
regular loose end, at which a blow-up yields the cracktip.
\end{abstract}
\maketitle

\section{Introduction}
\label{sec:introduction}

For a closed set $K\subset\R^2$ and a real-valued function
$u\in W^{1,2}_{\mathrm{loc}}(\R^2\setminus K)$, the Mumford--Shah energy
in a bounded open set $O\subset\R^2$ is
\[
 \E(u,K;O)=\int_{O\setminus K}|\nabla u|^2\,dx+\HH^1(K\cap O)\, ,
\]
where $\HH^1$ denotes one-dimensional Hausdorff measure.
We will call $(u,K)$ admissible if
$\E(u,K;B)<\infty$ for every bounded open set $B$ and, following the literature (cf. \cite{DLF}) 
an admissible pair is called an \emph{absolute global minimizer} if
$\E(u,K;B)\leq\E(v,J;B)$ for every bounded open $B$ and every admissible
competitor $(v,J)$ agreeing with $(u,K)$ outside a compact subset of $B$. The existence and regularity of minimizers for the functional (including a fidelity term) has a long history and it is not possible to do full justice to it in a short introduction. We include in the bibliography of this note some of the most relevant contributions (see \cite{DGA88,DGCL,DMMS,Bon,Dav96,AP,AFP97,Rig,BD,MSregularity,DL,Dav05,KLM,DLFhigher,DPF,AM,DLFG}), but we refer to the monograph \cite{DLF} for a much more extensive one and for an account of the history of the problem.

A well known conjecture of De Giorgi \cite{DG91} states that (up to obvious symmetries) the only global minimizer is given by a suitable pair $(K_*, u_*)$ where $K_*$ is a halfline, and $u_*$ an explicit harmonic function given below.  This pair is called {\em cracktip} in the literature and its global minimality was established 25 years ago
by Bonnet and David in their beautiful monograph \cite{BD}. The following is the precise statement.

\begin{theorem}[Bonnet--David]\label{thm:cracktip}
The pair
\[
 K_* =(-\infty,0]\times\{0\},\qquad
 u_*(r\cos\theta,r\sin\theta)
       =\sqrt{\frac2\pi}\,\sqrt r\sin\frac\theta2,
 \quad r>0,\quad -\pi<\theta<\pi,
\]
is an absolute global minimizer of the Mumford--Shah functional.
\end{theorem}

 It has been known since some time that the famous Mumford-Shah conjecture is in fact implied by the remaining part of De Giorgi's conjecture, namely by the absence of {\em any other} global minimizer (which in the literature have been called {\em exotic minimizers}). For an account of this, we refer to the monograph \cite{DLF}, which gives a complete description of the regularity theory developed in the literature for the planar problem since \cite{MS}. A proof of the nonexistence of exotic minimizers has been given recently in \cite{Da} using the assistance of chatgpt-astra-6 and this leads to a simple and elegant resolution of the last step towards the Mumford-Shah conjecture. The PhD thesis of the first author, \cite{DPhD}, whose original contributions are contained in the paper \cite{B24}, proposed a conjecture analogous to that of De Giorgi pertaining the Dirichlet boundary regularity theory. It turns out that this too can be proved combining insight gained through chatgpt-astra-6 and some further considerations due to the first author, cf. \cite{Db}. The overall result is that we have now a complete regularity theory for the Mumford-Shah functional both in the interior and at the boundary, as long as the latter is suitably smooth (\cite{B24,Db} covers the case of Dirichlet boundary conditions, while \cite{DL} covers the case of Neumann boundary conditions).

With the regularity theory is then easy to provide a very short and clean argument to show Theorem~\ref{thm:cracktip}. The key idea is to give a Dirichlet boundary condition on the unit disk $B=B_1 (0)$ which is smooth but forces the creation of a nontrivial set $K$ for the minimizer, because the harmonic competitor would have too large Dirichlet energy. The boundary and interior regularity theory would then imply that, if $K$ lacked an ``internal terminal point'', $B\setminus K$ would have a connected component which is ``shielded'' from the boundary, in other words an ``internal pocket''. In such  a pocket the minimizer can take an arbitrary constant value and this is easily seen to create a contradiction. But then if $K$ has an interior terminal point $p$, the interior regularity theory implies that its ``blow-up'' at $p$ must be at the same time a global minimizer and a cracktip. 

While the idea is simple, we want to underline that this does not really substitute the whole book \cite{BD}, which is full of useful and beautiful tools employed in the interior regularity theory (and which has inspired also the boundary regularity theory \cite{B24}).  In fact the first part of \cite{BD} is devoted to prove a partial ``uniqueness result'', which is employed in \cite{DLF} in order to prove the relevant $\varepsilon$-regularity theorem at the interior terminal point $p$. This does not require the minimality of the cracktip
(see the discussion preceding \cite[Theorem~1.4.6]{DLF}) and it is contained in, roughly, the first half of the book \cite{BD}.
The argument given in this note avoids therefore the second half of \cite{BD}, which can be thought as a way to ``force the right boundary conditions at planar infinity'' to show minimality of the cracktip. Since the cracktip has infinite energy on the entire plane, implementing the part just described requires a rather subtle modification of the Mumford-Shah functional and several annoying technicalities (which the authors of \cite{BD} resolve brilliantly but with a lot of work). The local argument given here is instead much more transparent. Clearly, compared to \cite{BD} there is a price to pay, which is the development of a suitable regularity theory for the Dirichlet problem. The latter is however a rather natural theory, which has its own independent interest. 

\section{Proof of Theorem~\ref{thm:cracktip}}
\label{sec:proof}

\subsection{Choice of a suitable boundary datum}
Let $B=B_1(0)$, choose any smooth function $\bar{g}$ on $\partial B_1 (0)$ and let $g(x)=2x_1$ be the boundary datum.
By the Dirichlet existence theorem \cite{CL}, we can fix a normalized
minimizer $(u,K)$ of
\[
 F(u,K)=\int_{B\setminus K}|\nabla u|^2\,dx+\HH^1(K),
 \qquad K\subset\overline B,\quad
 u\in W^{1,2}(B\setminus K),\quad u=g\text{ on }\partial B\setminus K.
\]
Here $K$ is closed and normalized means
$\HH^1(K\cap B_\rho(x))>0$ for every $x\in K$ and $\rho>0$.
By truncation we may assume $|u|\leq2$, cf. \cite[Section 2.1.1]{DLF}, while
the competitor $(0,\partial B)$ gives $F(u,K)\leq2\pi$.
If $K=\varnothing$, the harmonic function $g(x)=2x_1$ minimizes
the Dirichlet energy among functions with the prescribed boundary trace,
so $F(u,\varnothing)\geq\int_B|\nabla g|^2\,dx=4\pi$.
Thus $K\ne\varnothing$.

\subsection{Regularity theory}
By the resolution of the Mumford-Shah conjecture, \cite[Theorem~1.1]{Da}, every point of $K\cap B$
is a pure jump point, a regular loose end, or a triple junction,
in the sense of \cite[Definition~1.3.2]{DLF}.
Moreover the boundary classification in \cite[Theorem~2.1]{Db} and the
boundary $\varepsilon$-regularity \cite[Theorem~1.1]{B24} imply that,
near each point of $K\cap\partial B$, the set $K$ is the graph of a
$C^{1,\alpha}$ function over the tangent line to $\partial B$.
Consequently, no loose end or triple junction lies on $\partial B$.
By compactness, $K$ consists of finitely many $C^1$ arcs meeting at
triple junctions, together with possibly $C^1$ simple closed curves.

\subsection{Absence of pockets}
We adapt the argument given in
\cite[Lemma~4.3.1]{DLF} to exclude that $K$ can contain loops. Suppose indeed that any simple closed curve
$\Gamma$ (not necessarily $C^1$) is contained in $K$, and let $D$ be the bounded connected
component of $\R^2\setminus\Gamma$. Since $\Gamma\subset\overline B$,
we have $D\subset B$ and $\overline D\cap\partial B\subset\Gamma$.
Replacing $u$ by a constant on $D$ and $K$ by $K\setminus D$ is
admissible. By minimality and normalization, $K\cap D=\varnothing$
and $u$ is constant on $D$, with this constant freely adjustable
without affecting minimality.

If $\Gamma\cap B=\varnothing$, then $\Gamma=\partial B$ and $D=B$,
hence $K=\partial B$, contradicting \cite[Proposition~A.2]{B24}.
Thus $\Gamma\cap B\ne\varnothing$. Partial regularity gives a point
$x\in\Gamma\cap B$ where
$K$ is a regular arc and the exterior restriction of $u$ has a
$C^1$ extension $f$ across it
\cite[Corollary~3.1.5 and Lemma~4.3.1]{DLF}.

Replacing $u$ on $D$ by the constant $f(x)$ gives another minimizer.
The local comparison in the proof of \cite[Lemma~4.3.1]{DLF},
applied in a sufficiently small ball compactly contained in $B$,
then yields an admissible competitor with strictly smaller energy,
a contradiction.

\subsection{Blow-up at an interior terminal point}
Since $K$ is nonempty and contains no simple closed curve, it has
a loose end $p$: otherwise, following its finitely many arcs without
turning back would eventually revisit a point and produce a simple curve, as excluded in the previous paragraph.
The boundary regularity then implies that $p$ must be in the interior and the interior regularity implies that 
$K$ is locally a $C^1$ arc ending at $p$.
By blow-up compactness and Bonnet's monotonicity formula with its
equality case \cite[Theorem~2.2.3 and Proposition~2.6.1]{DLF},
a blow-up at $p$ is a generalized global minimizer
(in the sense of \cite[Definition~2.2.4]{DLF}), which (after a suitable rotation, addition of a constant, and a possible change of sign) is equal to the cracktip.

\section*{Acknowledgements}
The work of F. Deangelis was funded by the Deutsche Forschungsgemeinschaft
(DFG, German Research Foundation) under Germany's Excellence Strategy
EXC~2044/2 -- 390685587, Mathematics M\"unster:
Dynamics--Geometry--Structure.

\end{document}